\documentclass[preprint,review,12pt]{elsarticle}

\usepackage{graphics}
\usepackage{graphicx}
\usepackage{epsfig}

\usepackage{amssymb}

\begin{document}

\begin{frontmatter}



\title{Endpoints of  multi-valued weak  contractions on the  metric space valued in partially ordered groups}

\author{Congdian Cheng}

\address{College of Mathematics and Systems
Science, Shenyang Normal University,
 Shenyang,  110034,  P.R. China
{\bf Tel:}024-86170136; {\bf Fax:} 024-86171095 {\bf E-mail:}
 zhiyang918@163.com}


\begin{abstract}
We introduce the  metric space valued in partially ordered groups,
and define the convergence of sequences
 and
the multi-valued weak contractions, etc. , on the space. We then
establish  endpoint theorems for the defined maps. Our contributions
extend the theory of cone metric space constructed by Huang and
Zhang (2007) and some recent results on  the fixed point and
 endpoint theory, such as  the endpoint theorem given by Amini-Harandi
(2010).\end{abstract}


\end{frontmatter}


\section{Introduction}
\label{1}Let $(X, d)$ be a complete metric space.  Denote by $CB(X)$
the class of all nonempty closed and bounded subsets of $X$.  Denote
by $H(A,B)$ the Hausdorff metric of $A$ and $B$ with respect to $d$,
that is,
$$H(A,B)=\max\{\sup \limits_{x\in A}d(x,B), \sup \limits_{y\in B}d(y,A)\},$$
for all $A,B\in CB(X)$, where $d(x,B)=\inf\limits_{y\in B}d(x,y)$.
Further let $T : X\rightarrow 2^X$  be a multi-valued / set-valued
map. A point $x$ is called a fixed point of $T$ if $x\in Tx$. Define
$Fix(T) = \{x\in X : x \in Tx\}$.
A point $x$ is called an endpoint / a stationary point of a
multi-valued map $T$ if $Tx = \{x\}$. We denote the set of all
endpoints of $T$ by $End(T)$.

The investigation of endpoint of multi-valued mappings is an
important extending  of the study of  fixed point, which
was made as early as 30 years ago, and has received great attention
in recent years, see e.g. [1-2] and the references therein.
In particular, Amini-Harandi [1]
 (2010)  proved the Theorem 1.1 below.\\

 \noindent{\bf
Theorem 1.1} ([1, Theorem 2.1]).\textit{ Let $(X, d)$ be a complete
metric space and $T: X\rightarrow CB(X)$ be a set-valued map that
satisfies
$$H(Tx, Ty)
\leq \psi(d(x, y)),\eqno (1.1)$$ for each $x, y \in X$ , where $\psi
: [0,+\infty) \rightarrow [0,+\infty)$ is upper semicontinuous
(u.s.c.),  $\psi(t) < t$ for each $t
> 0$  and  satisfies $\liminf \limits_{
t\rightarrow +\infty }(t - \psi(t)) > 0$. Then $T$ has a unique
endpoint if and only if T has the approximate endpoint property. (
i.e. $\inf\limits_{x\in X}\sup\limits_{y\in Tx}d(x,y)=0$.
)}\\

Huang and Zhang [3]
(2007) introduced the concept of cone metric space, and  established
some fixed point theorems for contractive type maps in a normal cone
metric space. Subsequently, some other authors gave many results
about the fixed point theory in cone metric spaces. For example,
Rezapour and Hamlbarani [4] (2008)
 generalized some results of [3].
  Raja and Vaezpour [5]
(2008) presented some extensions of Banach$'$s Contraction Principle
in complete cone metric spaces. Aage and Salunke [6] (2011) proved
some fixed point theorems for the expansion onto mappings on
complete cone metric spaces. Also, many common fixed point theorems
were proved for maps on cone metric spaces in some literatures, for
example, see Ili\'{c} and Rako\v{c}evi\'{c} [7] (2008); Arshad, Azam
and Vetro [8] (2009), whose results generalized and unified many
fixed point theorems. Rezapour and Haghi [9] (2009), as well as
Haghi and Rezapour [10] (2010) studied fixed points of
multifunctions (i.e. multi-valued mappings) on normal cone metric
spaces and on regular cone metric spaces, respectively. Moreover,
Wardowski [11] (2009) introduced a kind of set-valued contractions
in cone metric spaces and established endpoint and fixed point
theorems for his contractions.

 In addition, Rezapour and Haghi [9]
(2009) introduced the concept of cone topology on cone metric space.
Lakshmikantham and \'{C}iri\'{c} [12] (2009) introduced the concept
of a mixed g-monotone mapping and prove coupled coincidence and
coupled common fixed point theorems for such nonlinear contractive
mappings in partially ordered complete metric spaces. Harjani and
Sadarangani [13] (2009) present some fixed point theorems for weakly
contractive maps in a complete metric space endowed with a partial
order. And Zhang [14] (2010)  proved some new fixed point and
coupled fixed point theorems for multivalued monotone mappings in
ordered metric spaces. Finally, Amini-Harandi [15] (2011) studied
fixed point theorems for a kind of generalized quasicontraction maps
in so called the vector modular spaces.

Motivated by the contributions stated above, the present work
introduces the  metric space valued in a partially ordered group
endowed with a topological structure and the metric space valued in
a partially ordered module endowed with a topological structure, and
establishes some fundamental concepts of analysis on the introduced
spaces, such as the convergence of sequences, which extends the
theory of cone metric space.  It also defines multi-valued weak
contractions, etc., on the introduced spaces. And then it focus on
addressing the endpoint theory of the multi-valued weak
contractions.


\section{Preliminaries}
\label{1} This section provides necessary preliminaries for our
 discussions.

We first make the following explanations. For a partial order
$\preceq$ of a set, we write $a\prec b$ to indicate that $a\preceq
b$ but $a\neq b$, where $a$ and $b$ are elements of the set. And for
a group $G$ with partial order $\preceq$, we write $G_+$ and $G^+$
to  indicate respectively the sets $\{a\in G:a\succeq \theta\}$ and
$\{a\in G:a\succ\theta\}$,
where $\theta$ indicates the identity element of $G$.\\

\noindent{\bf Definition 2.1.} Let $G$ be an abelian/a commutative
group with partial order $\preceq$. We call $G$  a
$\preceq$-partially ordered group, a partially ordered group for
simplicity,
 if $\prec$ satisfies the  law (g1)
$a\prec b \Rightarrow a+c\prec b+c, \forall a,b,c\in G$. Let further
$G$ be an $R$-module and the integral ring $R$ be a $\leq$-partially
ordered group.  Assume that the partial order $<$
satisfies  the  law (r1): $1>0$, where $1$ and $0$ are the unit
element and the identity element of $R$, respectively.
 Assume also that the partial orders $<$ and $\prec$ satisfy  the  law
 (m1): $a\prec b \Rightarrow ra\prec rb, \forall a,b\in G$ and
$\forall r\in R^+$ (i.e. $r>0$).
Then we call $G$  an $(R,\leq,\preceq)$-partially ordered module, a
partially ordered module for
simplicity.\\

\noindent{\bf Remark 2.2.} (i)For convenience, we focus our
attention to study under the assumption that there exist
non-identity elements in group $G$ below. (ii) Note that each
element of a group has an inverse element. From (g1),
we can easily obtain  the order relation: (g1)$'$ $a\preceq b
\Leftrightarrow a+c\preceq b+c, \forall a,b,c\in G$.
  In
addition, from (m1), we can easily obtain  the order relation:
(m1)$'$ $a\preceq b \Rightarrow ra\preceq rb, \forall a,b\in G$ and
$\forall r\in R_+$. (iii) From (m1), we can also obtain the order
relations: (m2)  $r< s\Rightarrow ra\prec sa, \forall r,s\in R$ and
$\forall a\in G^+$; and (m2)$'$  $r\leq s\Rightarrow ra\preceq sa,
\forall r,s\in R$ and $\forall a\in G_+$. In fact, let $r< s$, then,
by (g1), we have $0< s-r$. Let also $a\in G^+$, i.e. $\theta\prec
a$. Then, by (m1), we have $(s-r)\theta\prec (s-r)a\Rightarrow
\theta\prec (s-r)a$. From (g1), this leads to $ra\prec ra+(s-r)a=
sa$. So we have (m2). Finally, from (m2), it is obvious that we have
(m2)$'$.
 (iv) It is
obvious that the partially ordered
module is a special kind of the partially ordered group.\\

\noindent{\bf Example 2.3.} Let $E$ be a Banach space over the  real
field $\mathbb{R}$ and $P$ be a subset of $E$. $P$ is called a cone
if and only if: (i) $P$ is closed, nonempty, and $P\neq\{\theta\}$;
(ii) $ r, s \in{\mathbb{R}_+}, a,b \in P\Rightarrow ra + sb \in P $;
(iii) $ a \in P$ and $-a \in P \Rightarrow a = \theta$. Here
${\mathbb{R}}_+$ denotes all the non-negative real numbers. For a
given cone $P$ of $E$, define the partial order $\preceq$ on $E$ by
$x\preceq y $ if and only if $y-x\in P$, see [3]. Then it can be
easily verified that $E$ is an $(\mathbb{R},\leq,\preceq)$-partially
ordered module, and therefore, of course, is a $\preceq$-partially
ordered group. Here $\leq$ is the usual order of $\mathbb{R}$.\\

In the following part of this section $G$ is supposed either of a
$\preceq$-partially ordered group and
an $(R,\leq,\preceq)$-partially ordered module unless
otherwise specified.\\

\noindent{\bf Definition 2.4.}   Let $\ll$  be a non-empty relation
of $G$.
$\ll$ is called an analytic  topological structure of partially
ordered group $G$ if it satisfies:
 (t1) $a\ll b\Rightarrow a\prec b, \forall a, b
\in G$; (t2) $a\preceq b, b\ll c\Rightarrow a\ll c$; (t3)  $a\ll b
\Rightarrow a+c\ll b+c, \forall a,b,c\in G$; (t4) $\theta \preceq
a\ll \varepsilon, \forall \varepsilon\gg \theta \Rightarrow
a=\theta$; and (t5) $\forall \varepsilon\gg\theta$, there exists
$\eta\gg\theta$ such that $\eta\ll\varepsilon$. $\ll$ is called an
analytic topological structure of partially ordered module $G$ if it
also satisfies: (t6) $a\ll b \Rightarrow ra\ll rb, \forall
a,b\in G$ and $\forall r\in R^+$.\\

\noindent{\bf Remark 2.5.}
In the definition above, for $\ll$  is non-empty,  there are
actually infinite elements $\varepsilon$ such that
$\varepsilon\gg\theta$ in $G$. In fact, since $\ll$  is a non-empty
relation, there exist at least two elements $a$ and $b$ such that
$a\ll b$. By (t3),  we have $\theta\ll b-a$. Thus,
according to (t5), the result holds.  \\

\noindent{\bf Example 2.6.} For the partially ordered module $E$ of
Example 2.3, define the relation $\ll$ by $x\ll y$ if and only if
$y-x\in intP$, where $intP$ denotes the interior of $P$,
see [3] and [4].
Then we can  verify that $\ll$ is an analytic topological structure
of $E$. In fact, it is obvious that $\ll$ satisfies (t1), (t3), (t5)
and (t6).  To prove (t2), let $a\preceq b$ and $b\ll c$. Then, from
$b\ll c$, we have $\theta\ll (c-b)$. So there is an $r\in
  {\mathbb{R}}^+$ such that $N((c-b), r)\subset P$, where $N((c-b), r)
  =\{x\in E: \|x-(c-b)\|<r\}$ and $\|x\|$ indicates the norm of $x$. Consider $N((c-a),
  r)$.
Let $u\in N((c-a),
  r)$. Then $\|u-(b-a)-(c-b)\|=\|u-(c-a)\|<r$. This implies
  $u-(b-a)\in N((c-b), r)\subset
  P$. On the other hand, $(b-a)\in P$ for $a\preceq b$. So
  $u=u-(b-a)+(b-a)\in P$. Namely  $N((c-a),
  r)\subset P$. Hence $\theta\ll c-a$, e.g. $a\ll c$, that is, (t2)
holds. To prove (t4), assume $\theta \preceq a\ll \varepsilon,
\forall \varepsilon\gg \theta$. Let $c\gg\theta$. Then
$\frac{c}{n}\gg\theta$. By regarding $\frac{c}{n}$ as $\varepsilon$,
we have  $\frac{c}{n}\gg a\Rightarrow\frac{c}{n}-a\gg \theta
\Rightarrow\frac{c}{n}- a \in P$ for all $n\in \mathbb{N}$, where
$\mathbb{N}$ represents all the natural numbers. This leads to $- a
\in P$ because $\frac{c}{n}\rightarrow\theta$ (in norm) and $P$ is
closed. So, by $ a \in P$, we have $a=\theta$. That is, (t4) holds.
Therefore, $\ll$ is an analytic topological structure of $E$.\\

\noindent{\bf Definition 2.7.} Let $\ll$ be an analytic topological
structure of $G$ and $a\in G_{+}$. A sequence $\{a_{n}\}$ of $G_{+}$
is said to be convergent to $a$ (in $\ll$) if $\forall
\varepsilon\gg \theta$, there is a natural number $N$ such that
$\theta\preceq a_{n}-a\ll \varepsilon$ for all $n>N$, denoted
by $a_{n}\rightarrow a$ or $\lim\limits_{n\rightarrow\infty}a_n = a$.\\

\noindent{\bf Remark 2.8.} (i) Let $\ll$ be an analytic topological
structure of $G$, which is different from the $\prec$. Suppose
$\prec$ is also an analytic topological structure of $G$, and
sequence $\{a_{n}\}$ of $G_{+}$ converges to $\theta$ in $\prec$.
 Then we can easily know that $a_{n}$ converges to $\theta$ in $\ll$.
(In fact, let $\varepsilon\gg\theta$. Then, from (t5), there exists
$\eta\gg\theta$ such that $\eta\ll\varepsilon$. For the $\eta$,
since $a_{n}\rightarrow \theta$ in $\prec$, there is a natural
number $N$ such that $a_{n}\prec\eta$ for all $n>N$. By (t2) and
$\eta\ll\varepsilon$, this leads to $a_{n}\ll\varepsilon$ for all
$n>N$. Hence $a_{n}\rightarrow \theta$ in $\ll$.) That is, the
convergence in $\prec$ is stronger than in $\ll$. So, in the case,
the convergence in $\ll$ can be regarded as a kind of weak
convergence. (ii) For the analytic topological structure $\ll$ of
the partially ordered module $E$ in Example 2.6, it can be easily
verified that $\ll$ is different from the $\prec$ if $E$ is a
two-dimensional Euclidean space and $P=\{(x,y): x\geq 0, y\geq 0\}$.
(iii) It can be easily verified that for an analytic topological
structure $\ll$ of $G$, $a_{n}\rightarrow a$ $\Leftrightarrow$
$\forall b\gg a\succeq\theta$, there is a natural number $N$ such
that $a\preceq a_{n}\ll b$ for all $n>N$. In fact, assume
$a_{n}\rightarrow a$. $\forall b\gg a$, let $\varepsilon= b-a$.
Then, by (t3), we have $\varepsilon\gg\theta$. So, there is a
natural number $N$ such that $\theta\preceq a_{n}-a\ll \varepsilon$
for all $n>N$. From (g1)$'$ and (t3), this leads to $a\preceq
a_{n}\ll b$ for all $n>N$. Conversely,
$\forall\varepsilon\gg\theta$, let $ b=a+\varepsilon$. Then
 $ b\gg a\succeq\theta$. So, there is a natural number $N$ such that
$a\preceq a_{n}\ll b\Rightarrow\theta\preceq a_{n}-a\ll \varepsilon$
for all $n>N$. Note that $a\succeq\theta$.
This shows  $a_{n}\rightarrow a$.\\

\noindent{\bf Remark 2.9.} For  the $E$ and the  analytic
topological structure $\ll$ of Example 2.6, let $\{a_{n}\}$ be a
sequence in $E_+$. Assume
 $a_{n}\rightarrow \theta$ in norm. Then,
$\forall\varepsilon\gg\theta$, there exists $r\in \mathbb{R}^+$ such
that $N(\varepsilon,r)\subset P$. Due to  $a_{n}\rightarrow \theta$
in norm, there exists also a natural number $N$ such that
$\|a_{n}\|<r$ for all $n>N$. Therefore,
$\|(\varepsilon-a_{n})-\varepsilon\|<r\Rightarrow(\varepsilon-a_{n})\in
N(\varepsilon,r)\Rightarrow (\varepsilon-a_{n})\in intP$, that is,
$a_{n}\ll \varepsilon$, for all $n>N$. This implies that
$a_{n}\rightarrow \theta$ in $\ll$ if $a_{n}\rightarrow \theta$ in
norm.\\

$G$ always associates with an analytic topological structure $\ll$
and the  convergence of the sequences of $G_{+}$
  is in $\ll$ are assumed below.\\

\noindent{\bf Lemma 2.10.} \textit{Let $\{a_{n}\}$ and $\{b_{n}\}$
be two sequences of $G_{+}$. We have the three conclusions as
follows. (i) If $a_{n}\rightarrow\theta$, then
$\lim\limits_{n\rightarrow\infty}a_n $ is unique. (ii) If
$a_{n}\rightarrow\theta$ and $b_{n}\rightarrow\theta$, then
$a_{n}+b_{n}\rightarrow\theta$. (iii) If $b_{n}\succeq a_{n}
\succeq a\succeq\theta$ for all $n\in {\mathbb{N}}$ and
$b_{n}\rightarrow a$, then $(b_{n}-a_{n})\rightarrow\theta$.}\\

\noindent{\bf Proof}. Proving (i). Let
$\lim\limits_{n\rightarrow\infty}a_n =a$. Then there is a natural
number $N_{1}$ such that $a_{n}\succeq a$ for all $n>N_{1}$. On the
other hand, $\forall\varepsilon\gg\theta$, since
$a_{n}\rightarrow\theta$, there is a natural number $N_{2}$ such
that $\varepsilon\gg a_{n}\succeq\theta$ for all $n>N_{2}$. Let
$n=\max\{N_{1}, N_{2}\}+1$. Then, $\varepsilon\gg a_n$ and $a_n
\succeq a\succeq\theta$. From (t2), this leads to   $\varepsilon \gg
a\succeq\theta$.
 By virture of (t4), we
have $a=\theta$. Hence (i) holds.

Proving (ii).  Let $\varepsilon\gg\theta$. By (t5), there exists
$\eta\gg\theta$ such that $\varepsilon-\eta\gg\theta$. For
$a_{n}\rightarrow\theta$ and $b_{n}\rightarrow\theta$, there are
natural numbers $N_{1}$ and $N_{2}$ such that $a_n\ll \eta, \forall
n>N_{1}$ and $b_n\ll\varepsilon-\eta, \forall n>N_{2}$. Put
$N=\max\{N_{1},N_{2}\}$. We have: $a_n+b_n\ll
\eta+\varepsilon-\eta=\varepsilon, \forall n>N$. Hence,
$a_{n}+b_{n}\rightarrow\theta$. That is (ii) holds.

Proving (iii). Arguing by contradiction, assume
$\lim\limits_{n\rightarrow\infty}(b_{n}-a_{n})\neq\theta$.  Then
there exists a $\delta\gg\theta$ and a subsequence
$\{b_{n_i}-a_{n_i}\}$ such that $(b_{n_i}-a_{n_i})\ll\delta$ does
not hold for all $i\in {\mathbb{N}}$. From $a_{n_i}\succeq a$, by
(g1)$'$, we  have $b_{n_i}-a_{n_i}\preceq b_{n_i}-a$. This implies
that
 $b_{n_i}-a\ll\delta$ does not hold for all $i\in {\mathbb{N}}$. (In
 fact, if for some $i\in {\mathbb{N}}$,
 $b_{n_i}-a\ll\delta$, then, from  (t2) and  $b_{n_i}-a_{n_i}\preceq b_{n_i}-a$,
 we have $b_{n_i}-a_{n_i}\ll\delta$,
which contradicts that $(b_{n_i}-a_{n_i})\ll\delta$ does not hold.)
Hence $\lim\limits_{n\rightarrow\infty}b_{n}\neq a$. The
contradiction shows (iii) holds.
$\square$\\

\noindent{\bf Definition 2.11.}  $G$ is called regular if   every
decreasing sequence $\{a_n\}$ of $G_+$ is convergent. That is, if a
sequence $\{a_n\}$ of $G_+$  satisfies $a_n\preceq a_{n+1}$ for all
$n\in {\mathbb{N}}$, then
 exists a $a\in G_+$ such that
$a_n$ converges to $a$.\\

\noindent{\bf Remark 2.12.} (i)
Let $A\subseteq G_+, A\neq\emptyset$ and $a\in G_+$.  $a$ is called
the infimum of $A$ if and only if $a$ is a lower bound of $A$,  and
$c\preceq a$ for each lower bound $c$ of $A$, denoted by $inf A$.
(ii) It is obvious that there is
at most one infimum for each subset of $G_+$. In fact, for any
$A\subseteq G_+$, let $a$ and $b$ be two infimums of $A$. Then both
$a\preceq b$ and $b\preceq a$ hold. Hence $a=b$. This shows that $A$
has at most one infimum. (iii) In particular, $G$ is regular if for
each non-empty subset $A$ of $G_+$,
 $inf A$ exists and there exists a sequence $\{a_n\}$ of $A$ such that
$a_n$ converges to $inf A$. Actually, let $\{a_n\}$ be a decreasing sequence
 of $G_+$. Then, in the case,
   $inf \{a_n\}$ exists and there exists a subsequence $\{a_{n_i}\}$ of
  $\{a_n\}$ such that $\{a_{n_i}\}$ converges to $inf \{a_n\}$. Since $\{a_{n_i}\}$ converges to $inf \{a_n\}$,
  $\forall \varepsilon\gg\theta$, there is a natural number $I$
   such that $\varepsilon\gg a_{n_i}-inf \{a_n\}\succeq\theta$ for all $i\geq I$.
  For $\{a_n\}$ decreasing, this leads to
  $\varepsilon\gg a_{n}-inf \{a_n\}\succeq\theta$
 for all $n\geq n_I$. That is,  $\{a_{n}\}$ converges
 to $inf \{a_n\}$. Hence $G$ is regular.\\

\noindent{\bf Definition 2.13.} Let X be a non-empty
set. Suppose the mapping $d : X \times X \rightarrow G$ satisfies\\
(d1) $d(x, y)\succeq\theta$ for all $x, y \in X$ and $d(x, y) =
\theta$ if and
only if $x = y$,\\ (d2) $d(x, y) =d(y, x)$ for all $x, y \in X$,\\
(d3) $d(x, y)\preceq d(x, z)+d(z, y)$ for all $x, y, z \in X$.\\
Then d is called a metric ( on $X$) valued in partially ordered
group $G$, and $(X,d)$ is called a metric space valued in partially
ordered group $G$, when $G$ is a partially ordered group; a metric
of group and a metric space of group for simplicity, respectively.
(Then d is called a metric valued in partially ordered module $G$,
and $(X,d)$ is called a metric space valued in partially ordered
module $G$, when $G$ is a
partially ordered module.)\\

In the rest of this section we always assume that $(X,d)$ is either
of a metric space valued in partially ordered group $G$ and a metric
space
valued in partially ordered module $G$.\\

\noindent{\bf Definition 2.14.} For given $(X,d)$, let $x \in X$ and
$\{x_n\}$ be a sequence in $X$.
\\(i) We call that $\{x_n\}$ converges to $x$ if and only if
 $d(x_n,x)\rightarrow \theta$, denoted by
$\lim\limits_{n\rightarrow\infty}x_n = x$ or $x_n\rightarrow x$.
\\(ii) $\{x_n\}$ is a Cauchy sequence if and only if $d(x_n,
x_m)\rightarrow \theta$, that is, $\forall \varepsilon \gg\theta$,
 there is a natural number $N$ such that $d(x_n,
x_m)\ll\varepsilon$ for all $n,m \geq N$. \\(iii) $(X,d)$ is
complete
 if and only if every Cauchy
sequence is convergent. \\(iv) $(X,d)$ is regular
 if and only if $G$ is regular.\\

\noindent{\bf Remark 2.15.}
The relation between the regular space and the
complete space is an  interesting question for further research.\\

\noindent{\bf Definition 2.16.} Given $(X,d)$, let $T : X
\rightarrow (2^X-\emptyset)$ be
 a multi-valued mapping and $\varphi: X
\times X \rightarrow G_+$ be a  mapping with $\varphi(x,y)\prec
d(x,y)$ for all $d(x,y)\succ\theta$. $T$ is called  a multi-valued
($\varphi$-)weak contraction on $(X,d)$ if, for all different
$x,y\in X$, $\forall x'\in Tx$, there exists $\bar{y}\in Ty$ such
that
$$d(x',\bar{y})\preceq \varphi (x,y).\eqno (2.1)$$
$T$ is called a global multi-valued ($\varphi$-)weak contraction on
$(X,d)$ if, for all different $x,y\in X$,  we have
$$d(x',y')\preceq \varphi(x,y), \forall x'\in Tx, \forall y'\in
Ty.\eqno (2.2)$$ The weak contraction $T$ is called to satisfy
C-condition (convergence condition) if $d(x_n, y_n) - \varphi(x_n,
y_n)\rightarrow\theta$,
then $d(x_n, y_n)\rightarrow\theta$, where $x_n$ and $y_n$ are two sequences of $X$.\\

\noindent{\bf Remark 2.17.}  (i) It is obvious that a global
multi-valued weak contraction is a  multi-valued weak contraction.
(ii) The weak contraction $T$ is called to satisfy C$'$-condition if
$d(x_n, x_m) - \varphi(x_n, x_m)\rightarrow\theta (n\neq m)$,
 that is, $\forall \varepsilon\gg\theta$,
there is a $N$ such that $d(x_n, x_m)-\varphi(x_n,
x_m)\ll\varepsilon$ for all $n,m>N$ and $n\neq m$, then $d(x_n,
x_m)\rightarrow\theta$. (iii) If $T$ satisfies C-condition, then it
also satisfies C$'$-condition. In fact, for the set $\{(n,m): n,
m\in \mathbb{N}, n\neq m\}$ is countable, it can be rewritten as the
sequence $\{(x'_{i}, y'_{i})\}$. Assume $d(x_n, x_m) - \varphi(x_n,
x_m)\rightarrow\theta (n\neq m)$. Let $\varepsilon\gg\theta$. Then
there exists a natural number $N$ such that $d(x_{n},
x_{m})-\varphi(x_{n}, x_{m})\ll\varepsilon$ whenever $n,m>N$ and
$n\neq m$. Because the set $\{(n,m): n, m\leq N\}$ is finite, there
is a natural number $I$ such that if $i>I$ and $(x'_{i},
y'_{i})=(x_{n}, x_{m})$, then $n, m> N$. This implies $\forall i>I$,
we have $d(x'_{i}, y'_{i})-\varphi(x'_{i}, y'_{i})\ll\varepsilon$.
Hence  $d(x'_{i}, y'_{i})-\varphi(x'_{i}, y'_{i})\rightarrow\theta$.
For $T$ satisfies C-condition, we have $d(x'_{i},
y'_{i})\rightarrow\theta$. Further, due to $d(x'_{i},
y'_{i})\rightarrow\theta$, $\forall \varepsilon\gg\theta$, there is
 a natural number $I'$ such that $d(x'_{i},
y'_{i})\ll\varepsilon$ whenever $i>I'$. Since the set $\{i: i\leq
I'\}$ is finite, there is a natural number $N'$ such that if
$n,m>N'$, $n\neq m$ and $(x_{n}, x_{m})=(x'_{i}, y'_{i})$, then
$i>I'$. That is,  $\forall n,m>N'$, $n\neq m$, we have $d(x_n,
x_m)\ll\varepsilon$. Note that $d(x_n, x_m)=\theta$ when $n= m$.
This leads to $d(x_n,
x_m)\rightarrow\theta$.\\

\noindent{\bf Definition 2.18.} A map $T : X \rightarrow
(2^X-\emptyset)$ on $(X,d)$ is said to have approximate endpoint
property if there exist a sequence $\{x_n\}$ of $X$ and a sequence
$\{a_n\}$ of $G_+$ with $a_n\rightarrow\theta$ such that
$$d(x_n,x'_n)\preceq a_n, \forall x'_n\in Tx_n, \eqno (2.3)$$
for all $n\in {\mathbb{N}}$.\\

\noindent{\bf Remark 2.19.} When $(X,d)$ is the usual complete
metric space, it can be easily verified that $T$ has the approximate
endpoint property in Theorem 1.1 and $T$ has the approximate
endpoint property  defined in Definition 2.18 are equivalent. (In
fact, if $T$ has the approximate endpoint property in Theorem 1.1,
that is, $\inf\limits_{x\in X}\sup\limits_{y\in Tx}d(x,y)=0$, then
there is a sequence $\{x_n\}$ of $X$ such that
 $\sup\limits_{y\in Tx_n}d(x_n,y)\rightarrow0$. Let
  $a_n=\sup\limits_{y\in Tx_n}d(x_n,y)$. Then   $d(x_n,x'_n)\leq a_n, \forall x'_n\in Tx_n$ and
   $a_n\rightarrow0$. This shows that $T$ has the approximate
endpoint property  defined in Definition 2.18. On the other hand, if
$T$ has the approximate endpoint property  defined in Definition
2.18, that is, there exist a sequence $\{x_n\}$ of $X$ and a
sequence $\{a_n\}$ of ${\mathbb{R}}_+$ with $a_n\rightarrow 0$ such
that $d(x_n,x'_n)\leq a_n, \forall x'_n\in Tx_n$, then
$\sup\limits_{y\in Tx_n}d(x_n,y)\leq a_n$
 for all $n\in {\mathbb{N}}$. This implies that $\inf\limits_{x\in X}\sup\limits_{y\in Tx}d(x,y)=0$, namely,
 $T$ has the approximate endpoint property in Theorem 1.1. Hence we have the equivalence stated above.) \\

 \noindent{\bf Lemma 2.20.} \textit{
 Let $T$ be a multi-valued weak
  contraction on $(X,d)$. Then we have the following two conclusions.
  (i) $T$ has approximate
 endpoint property if $T$ has endpoints. (ii) $T$ has one endpoint at most.
 (i.e. $|End(T)|\leq 1$. Here $|End(T)|$
denotes the cardinal number of $End(T)$.)}\\

\noindent{\bf Proof}.  (i) is obvious. In fact, let $x$ be an
endpoint of $T$. Put $x_n=x$ and $a_n=\theta$ for all $n\in
{\mathbb{N}}$. Then $a_n\rightarrow\theta$ and (2.3) holds for all
$n\in {\mathbb{N}}$. Hence $T$ has the approximate endpoint
property. To prove (ii), assume $|End(T)|>1$. Then, there exist
$x,y\in End(T)$ such that $x\neq y$. From (2.1), we have
$d(x,y)\preceq\varphi(x,y).$
Note that $\varphi(x,y)\prec d(x,y)$ for any $d(x,y)\succ\theta$.
This implies  $d(x,y)=\theta$. Hence, from (d1), we have $x=y$. This
contradicts
$x\neq y$. So $|End(T)|\leq 1$, that is, (ii) holds. $\square$\\


\section{Main results}

In this section, we always assume that $(X,d)$ is a metric space
valued in partially ordered group $G$.

Now we are ready to  prove our main results. We first present the
following Theorem 3.1,  which extends Theorem 1.1 (Amini-Harandi [1,
Theorem 2.1]) to the
case of the metric space of group.\\

\noindent{\bf Theorem 3.1.} \textit{Let
 $T$ be a multi-valued weak
contraction on complete $(X,d)$ and
satisfy   C-condition. Then $T$ has a unique endpoint if and only if
it has the
approximate endpoint property.}\\

\noindent{\bf Proof.} The necessity is clear from the (i) of Lemma
2.20. Next we prove the sufficiency.

Since $T$ has the approximate endpoint
 property, there exist sequences $\{x_n\}$ of $X$ and $\{a_n\}$ of $G_+$ satisfying
 (2.3) and $a_n\rightarrow \theta$. If there exists a subsequence $\{x_{n_i}\}$
 such that $x_{n_i}$ being the same point $x$ of $X$ for all $i\in{\mathbb{N}}$, then we can easily know
 that $x$ is an endpoint of $T$ from (2.3). (In fact, for any given $x'\in
 Tx$, we have $d(x, x')\preceq a_{n_i}$ for all $i\in {\mathbb{N}}$. Since  $a_n\rightarrow
 \theta$, we have $a_{n_i}\rightarrow
 \theta(i\rightarrow\infty)$. So we have $\theta\preceq d(x, x')\ll \varepsilon$ for all
  $\varepsilon\gg\theta$. This implies $d(x, x')=\theta$ from (t4), that is,  $x=x'$.  Hence $x$ is an endpoint of
 $T$.)
 Otherwise, without loss of generality, we can assume $x_n\neq x_m$ whenever
 $n\neq m$ and continue to prove as follows.


For any different $n,m\in\mathbb{N}$, let $x'_n\in Tx_n$. Then
$$d(x_n,x_m)\preceq d(x_n,x'_n)+d(x'_n,x_m). \eqno(3.1)$$
Since $x_n\neq x_m$, according to  (2.1), there exists $\bar{x}_m\in
Tx_m$ such that
$$d(x'_n,\bar{x}_m)\preceq \varphi(x_n,x_m).$$
Using this and (3.1), we further obtain
$$\begin{array}{rcl}d(x_n,x_m)&\preceq & d(x_n,x'_n)+d(x'_n,\bar{x}_m)+d(\bar{x}_m,x_m)
\\&\preceq& d(x_n,x'_n)+\varphi(x_n,x_m)+d(\bar{x}_m,x_m).\end{array}\eqno(3.2)$$
By (2.3), we have $d(x_n,x'_n)\preceq a_n$ and
$d(\bar{x}_m,x_m)\preceq a_m$. From (3.2), this leads to
$$\begin{array}{rcl}& &d(x_n,x_m)\preceq a_n+\varphi(x_n,x_m)+a_m\\&\Rightarrow&
d(x_n,x_m)-\varphi(x_n,x_m)\preceq a_n+a_m, n\neq m.
\end{array}\eqno(3.3)$$
For $x_n\neq x_m$, we have $d(x_n,x_m)\succ\theta$. So,
$d(x_n,x_m)-\varphi(x_n,x_m)\succeq\theta$. On the other hand,
noting that $a_n\rightarrow\theta$, following the proof on the (ii)
of Lemma 2.10, we can easily know that $a_n+a_m\rightarrow \theta$.
Thus we can obtain $d(x_n,x_m)-\varphi(x_n,x_m)\rightarrow \theta
(x_n\neq x_m)$ from (3.3). This implies $d(x_n,x_m)\rightarrow
\theta$ for $T$ satisfies the C-condition, which leads to $\varphi$
satisfies the C$'$-condition, see the (ii) and (iii) of Remark 2.17.
Hence, $\{x_n\}$ is a Cauchy sequence. Since $(X, d)$ is  complete,
there is a $x\in X$ such that $x_n\rightarrow x$.

We show $x$ is an endpoint of $T$ below.

Since $x_n\neq x_m$ whenever $n\neq m$, without loss of generality,
we can assume $x_n\neq x$ for any $n\in{\mathbb{N}}$. Let $x'\in
Tx$. Then, for all $n\in{\mathbb{N}}$, we have
$$d(x',x)\preceq d(x',x_n)+d(x_n, x). \eqno(3.4)$$
For $x_n\neq x$, by (2.1), there exists $\bar{x}_n\in Tx_n$ such
that
$$d(x',\bar{x}_n)\preceq \varphi(x,x_n)\prec d(x,x_n). \eqno(3.5)$$
In terms of  (3.4), (3.5)  and  (2.3),  we obtain
$$\begin{array}{rcl}d(x',x)&\preceq& d(x',\bar{x}_n)+d(\bar{x}_n,x_n)+d(x_n,x)\\
&\preceq&\varphi(x,x_n)+a_n+d(x, x_n)\\&\preceq&d(x, x_n)+a_n+d(x,
x_n).\end{array}$$ Since $d(x, x_n)\rightarrow\theta$ and
$a_n\rightarrow\theta$, by the (ii) of Lemma 2.10, we obtain $d(x,
x_n)+d(x, x_n)+a_n\rightarrow \theta$. So,  $d(x,x')\ll\varepsilon$
for all $\varepsilon\gg\theta$. Note also that
$d(x,x')\succeq\theta$. From (t4), we have $d(x,x')=\theta$. Hence
$x=x'$. That is, $x\in End(T)$.

Finally, the uniqueness is directly obtained from the (ii) of Lemma
2.20. The
proof completes. $\square$\\

\noindent{\bf Remark 3.2.} Here we make a simple explanation for
Theorem 3.1  extending Theorem 1.1. Firstly, it is obvious that for
the usual order $\leq$ of the real field $\mathbb{R}$, $\mathbb{R}$
is a
$\leq$-partially ordered group with analytic topological structure
$>$. Further, due that
$(X,d)$ in Theorem 1.1 is a complete metric space, it is a complete
metric space of the group $\mathbb{R}$ with analytic topological
structure $>$.
 That is, $(X,d)$ satisfies
the requirement of Theorem 3.1. Secondly, for the $\psi(t)$ in
Theorem 1.1, let $\varphi(x,y)=\psi(d(x,y))$, then $\varphi(x,y)$ is
a mapping from $X \times X$ to $\mathbb{R}_+$ and
$\varphi(x,y)<d(x,y)$ for all $d(x,y)>\theta$.
For the mapping $T$
in Theorem 1.1 and the $\varphi$ defined above, we have $H(Tx,
Ty)\leq\varphi(x,y)$ for all different $x, y\in X$, that is,
$$\max\{\sup \limits_{x'\in Tx}d(x',Ty), \sup \limits_{y'\in
Ty}d(Tx,y')\}\leq\varphi(x,y).$$ This leads to $\sup \limits_{x'\in
Tx}d(x',Ty)\leq\varphi(x,y)\Rightarrow d(x',Ty)\leq \varphi(x,y),
\forall x'\in Tx$. Thus, for $Ty$ is  closed and bounded, $\forall
x'\in Tx$, there is a $\overline{y}\in Ty$ such that
$d(x',\overline{y})\leq \varphi(x,y)$.
This shows that
$T$  is  a multi-valued weak contraction on the space $(X,d)$.
Thirdly, if $d(x_n, y_n)-\varphi(x_n, y_n)\rightarrow0$, then
$d(x_n, y_n)\rightarrow 0$, that is, $T$ satisfies the C-condition.
In fact, if $d(x_n, y_n)$ does not converge to $0$, then there exist
a $\delta>0$ and a subsequence $\{d(x_{n_i}, y_{n_i})\}$ such that
$d(x_{n_i}, y_{n_i})-\psi(d(x_{n_i}, y_{n_i}))>\delta$ for all
$i\in{\mathbb{N}}$. $\lceil$ We show the fact is true as follows.
Let $d_n=d(x_n, y_n)$ and $d_n$ do not converge to $0$. If $\{d_n\}$
is unbounded, without loss of generality, we can assume that
$\{d_n\}$ increases and converges to $+\infty$. For $\liminf
\limits_{ t\rightarrow +\infty }(t - \psi(t))=\lim \limits_{
t\rightarrow +\infty }[\inf\{(s - \psi(s)): s>t\}]> 0$ and
$\inf\{(d_k - \psi(d_k)): k>n\}\geq\inf\{(s - \psi(s)): s>d_n\}$, we
have $\liminf \limits_{ n\rightarrow \infty }(d_n - \psi(d_n))=\lim
\limits_{ n\rightarrow \infty }[\inf\{(d_k - \psi(d_k)):
k>n\}]\geq\lim \limits_{ n\rightarrow \infty }[\inf\{(s - \psi(s)):
s>d_n\}]=\lim \limits_{ t\rightarrow +\infty }[\inf\{(s - \psi(s)):
s>t\}]> 0$.
Hence there exist a $\delta>0$ and a subsequence $\{d_{n_i}\}$ such
that $d_{n_i}-\psi(d_{n_i})>\delta$  for all $i\in{\mathbb{N}}$,
that is, the fact is true.
If $\{d_n\}$ is bounded, without loss of generality,  we  assume
that $\{d_n\}$ increases and converges to a point $t'>0$. Then, for
$\psi$ is u.s.c. at $t'$, i.e. $\limsup\limits_{ t\rightarrow t'
}\psi(t)\leq\psi(t')$,
 and $\psi(t') < t'$, we have $\liminf\limits_{ t\rightarrow t'
}(t-\psi(t))=t'-\limsup\limits_{ t\rightarrow t' }\psi(t)\geq
t'-\psi(t')>0$. Note that $\{d_n\}$ increases and $\lim
\limits_{\Delta t\rightarrow 0}[\inf\{(t - \psi(t)): 0<|t-t'|<\Delta
t\}]=\liminf\limits_{ t\rightarrow t' }(t-\psi(t))$. We have
$\liminf \limits_{ n\rightarrow \infty }(d_n - \psi(d_n))=\lim
\limits_{ n\rightarrow \infty }[\inf\{(d_k - \psi(d_k)):
k>n\}]\geq\lim \limits_{ n\rightarrow \infty }[\inf\{(t - \psi(t)):
d_n<t<t'\}]=\lim \limits_{\Delta t\rightarrow 0}[\inf\{(t -
\psi(t)): 0<t-t'<\Delta t\}]\geq\lim \limits_{\Delta t\rightarrow
0}[\inf\{(t - \psi(t)): 0<|t-t'|<\Delta t\}]=\liminf\limits_{
t\rightarrow t' }(t-\psi(t))>0$. Hence the fact is also
true.$\rceil$
 That is, $d(x_{n_i}, y_{n_i})-\varphi(x_{n_i}, y_{n_i})>\delta$ for all $i\in{\mathbb{N}}$.
This contradicts $d(x_n, y_n)-\varphi(x_n, y_n)$ converges to $0$.
Hence $T$ satisfies the C-condition.
Finally, for the $T$ has the approximate endpoint property of
Theorem 1.1, from Remark 2.19,
it has the approximate endpoint property (defined in Definition
2.18).
Hence we can directly obtain
Theorem 1.1  from Theorem 3.1.\\

Next we further present the following Theorem 3.3, which shows, in
the setting that $(X,d)$ is  complete and  regular, if the global
multi-valued weak contraction satisfies C-condition, then it
has the approximate endpoint property, so has a unique endpoint from Theorem 3.1.\\

\noindent{\bf Theorem 3.3.} \textit{Let $(X,d)$ be complete and
regular, $T$ be a global multi-valued weak contraction on $(X,d)$
and satisfy  C-condition.
Then $T$ has a unique endpoint.}\\

\noindent{\bf Proof}. We first prove the existence of  endpoints.

Arguing by contradiction, assume $T$ has no endpoint. Then for any
$x\in X$, there is at least one $y\in Tx$ such that $y\neq x$. Hence
 there must be a sequence $\{y_n\}$ of $X$ such that
$y_{n+1}\in Ty_n$ and $y_{n+1}\neq y_n$ for all $n\in \mathbb{N}$.
Note $T$ is a global multi-valued weak contraction.  In terms of
$y_{n+1}\in Ty_n$, $y_{n+1}\neq y_n$, (2.2) and $\varphi(x,y)\prec
d(x,y)$ for $d(x,y)\succ \theta$, we have
$$d(y_{n+1}, y_{n+2})\preceq \varphi(y_{n}, y_{n+1})\prec d(y_{n}, y_{n+1})\eqno(3.6)$$
for all $n\in\mathbb{N}$. Hence the sequence
  $\{d(y_{n}, y_{n+1})\}$ is decreasing. So, for $G$ is
regular, there exists $a\in G_+$ such that $d(y_{n},
y_{n+1})\rightarrow a$. Hence, from (3.6), we have $a\preceq
\varphi(y_{n}, y_{n+1})\prec d(y_{n}, y_{n+1})$. Further, according
to the (iii) of Lemma 2.10, we  obtain $d(y_{n},
y_{n+1})-\varphi(y_{n}, y_{n+1})\rightarrow \theta$.  For $T$
satisfies  C-condition, this leads to $d(y_{n}, y_{n+1})\rightarrow
\theta$. And by (3.6) we further have  $\varphi(y_{n},
y_{n+1})\rightarrow \theta$. Now let $x_{n}=y_{n+1}$ and
$a_{n}=\varphi(y_{n}, y_{n+1})$. Then we have $d(x_{n},
x'_{n})\preceq a_{n}, \forall x'_{n}\in Tx_{n}$ and
$a_{n}\rightarrow \theta$. That is, $T$ has the approximate endpoint
property. Thus, by Theorem 3.1, $T$ has endpoints. This contradicts
our assumption.  Hence the existence of endpoints is true.

Finally, the uniqueness follows directly from the (ii) of Lemma
2.20. This ends the proof.
$\square$\\

For the single-valued weak contraction can be regarded as a kind of
specific global multi-valued weak contraction, from Theorem 3.3, we
can immediately derive the Corollary 3.4 below, which generalizes
Lemma 2.4 and Corollary 2.5 of [1].\\

\noindent{\bf Corollary 3.4.} \textit{Let $(X,d)$ be complete and
regular,
 $f$ be a single-valued weak
contraction on $(X,d)$ and satisfy  C-condition. Then $T$ has a
unique fixed point.}\\

\noindent{\bf Remark 3.5.}  The multi-valued weak contraction can
not have the approximate endpoint property, even in the usual metric
space, for instance, see the Example 2.3 of [1].
\section{Endpoint theory for the metric space of module}

Note that a metric space of module is  a special metric space of
group. As applications of the  results proved above,  this section
 discusses the endpoint theory for the metric space of module. We
always assume  $(X,d)$ is a metric space of module
$G$ in this section.

By regarding $\alpha(x,y)d(x,y)$ as $\varphi(x,y)$,
we can easily derive the next Theorem 4.1 from Theorem 3.1 and Theorem 3.3.\\

\noindent{\bf Theorem 4.1.} \textit{Let $(X,d)$ be complete, $T : X
\rightarrow (2^X-\emptyset)$ be a multi-valued mapping. Let also
$\alpha: X\times X\rightarrow [0, 1) (=\{r\in R: 0\leq r< 1\})$
 be a  mapping, which satisfies:
 for any two
  sequences $\{x_n\}$ and $\{y_n\}$ of $X$ with $[1-\alpha(x_n,y_n)]d(x_n,y_n)\rightarrow\theta$,
  there is an $\alpha\in [0, 1)$
 such that $\alpha(x_n,y_n)\leq \alpha$ for all $n\in{\mathbb{N}}$,
 as well as
  $(1-\alpha)$ has multiplicative inverse
 $(1-\alpha)^{-1}$  and $(1-\alpha)^{-1}>0$.
 Then we have the following two conclusions.\\ (i) Suppose  for
all different $x,y\in X$, $\forall x'\in Tx$, there exists
$\bar{y}\in
 Ty$ such that
$d(x',\bar{y})\preceq \alpha(x,y)d(x,y)$.
 Then $T$ has a
unique endpoint if and only if $T$ has the approximate endpoint
 property.\\
(ii) Let $(X,d)$ be  regular. Suppose for all different $x,y\in X$,
we have $d(x', y')\preceq \alpha(x,y)d(x,y)$, $\forall x'\in Tx,
\forall y'\in Ty$.  Then $T$ has a
unique endpoint.}\\


\noindent{\bf Proof.}
Let  $\varphi(x,y)=
\alpha(x,y)d(x,y)$. Then
$\varphi: X\times X \rightarrow G_+$ is a mapping. Since
$\alpha(x,y)\prec 1$, by (m2), see the (iii) of remark 2.2, we have
 $\varphi(x,y)=\alpha(x,y)d(x,y)\prec d(x,y)$ for all  $d(x,y)\succ\theta$. We  prove the statement that
 if
$d(x_n, y_n)-\varphi(x_n, y_n)\rightarrow\theta$, then
$d(x_n,y_n)\rightarrow\theta$ below.

Let $\{x_n\}$ and $\{y_n\}$ be two sequences of $X$. Then $d(x_n,
y_n)-\varphi(x_n, y_n)=[1-\alpha(x_n,y_n)]d(x_n,y_n)$. So, if
$d(x_n, y_n)-\varphi(x_n, y_n)\rightarrow\theta$, then
$$[1-\alpha(x_n,y_n)]d(x_n,y_n)\rightarrow\theta. \eqno(4.1)$$
Hence there exists an $\alpha\in [0, 1)$   such that
$\alpha(x_n,y_n)\leq\alpha, \forall n\in{\mathbb{N}}$. By (g1)$'$,
this further leads to $(1-\alpha)\leq [1-\alpha(x_n,y_n)]$. Since
also $d(x_n,y_n)\succeq\theta$, by (m2)$'$, we have
$$(1-\alpha)d(x_n,y_n)\preceq[1-\alpha(x_n,y_n)]d(x_n,y_n), \forall n\in{\mathbb{N}}. \eqno(4.2)$$
In terms of  (4.1) and (4.2), we obtain
$(1-\alpha)d(x_n,y_n)\rightarrow\theta$. For $\alpha<1$, we have
$1-\alpha>0$. Let $\varepsilon\gg\theta$. By (t6), we have
$(1-\alpha)\varepsilon\gg\theta$. Hence,
 there exists a natural $N$ such that
$(1-\alpha)d(x_n,y_n)\ll(1-\alpha)\varepsilon$ for all $n>N$. For
$(1-\alpha)^{-1}>0$, from (t6), we obtain also
$d(x_n,y_n)\ll\varepsilon$ for all $n>N$. This implies
$d(x_n,y_n)\rightarrow\theta$.

For conclusion (i), by $\varphi(x,y)=\alpha(x,y)d(x,y)\prec d(x,y)$
for all $d(x,y)\succ\theta$ and the statement proved above, it is
obvious that $T$ is a multi-valued weak contraction on the space
$(X,d)$ and satisfies  C-condition. Hence we can immediately know
that the conclusion is true from Theorem 3.1. For conclusion (ii),
$T$ is clearly a global multi-valued weak contraction and satisfies
C-condition. Note that $(X,d)$ is regular. We can immediately know
that the conclusion is true from Theorem 3.3.
This completes the proof. $\square$\\

Replacing, in  Theorem 4.1, $\alpha(x,y)$ by $\alpha$, we
 directly obtain the following Corollary 4.2.\\

\noindent{\bf Corollary  4.2.} \textit{Let $(X,d)$ be complete,  $T
: X \rightarrow (2^X-\emptyset)$ be a multi-valued mapping. Let also
$\alpha\in [0, 1)$ with $(1-\alpha)^{-1}$ and $(1-\alpha)^{-1}>0$.
(i) Suppose $T$ satisfies for all different $x,y\in X$, $\forall
x'\in Tx$, there exists $\bar{y}\in
 Ty$ such that
$d(x',\bar{y})\preceq \alpha d(x,y)$. Then $T$ has a unique endpoint
if and only if $T$ has the approximate endpoint
 property. (ii) Let  $(X,d)$ be  regular. Suppose $T$ satisfies for all
  different $x,y\in X$,
$d(x',y')\preceq \alpha d(x,y), \forall x'\in Tx, \forall y'\in Ty$.
Then $T$ has a unique endpoint.}\\

\noindent{\bf Proof}. Let $\alpha(x,y)=\alpha$. And then applying
Theorem 4.1, we obtain  the Corollary instantly. $\square$\\

Finally,  in the (ii) of Corollary 4.2,  replacing also multi-valued
mapping by single-valued mapping, we obtain Corollary 4.3 below.\\

\noindent{\bf Corollary  4.3.} \textit{Let $(X,d)$ be complete and
 regular, $T : X \rightarrow X$ be a single-valued mapping,
$\alpha\in [0, 1)$ with $(1-\alpha)^{-1}$ and $(1-\alpha)^{-1}>0$.
Suppose $T$ satisfies for all different $x,y\in X$, $d(Tx,
Ty)\preceq \alpha d(x,y)$. Then $T$ has a unique fixed point.}\\

\noindent{\bf Remark 4.4.} In particular,   when $(X,d)$ is the
usual complete metric space, Corollary 4.3 is just the famous Banach
fixed point theorem.


\noindent{\bf Acknowledgements}

The author cordially thanks the anonymous referees for their
valuable comments which lead to the improvement of this paper.





\bibliographystyle{model1-num-names}
\bibliography{<your-bib-database>}







\end{document}